\documentclass[reqno]{amsart}
\usepackage{amssymb,stmaryrd}
\usepackage{amsfonts}
\usepackage{amstext}
\usepackage{algorithmic}
\usepackage{algorithm}
\usepackage{graphicx}
\usepackage{epstopdf}
\usepackage[all]{xy}
\usepackage{MnSymbol,slashed}

\numberwithin{equation}{section}
\newtheorem{theorem}{Theorem}[section]
\newtheorem{lemma}[theorem]{Lemma}
\newtheorem{proposition}[theorem]{Proposition}
\newtheorem{corollary}[theorem]{Corollary}

\theoremstyle{definition}

\newtheorem{example}[theorem]{Example}

\theoremstyle{remark}

\newcommand{\R}{{\mathbb{R}}}

\newcommand{\C}{{\mathbb{C}}}
\newcommand{\Z}{{\mathbb{Z}}}

\newcommand{\N}{{\mathbb{N}}}

\newcommand{\<}{{\langle}}
\renewcommand{\(}{{(}}

\renewcommand{\>}{{\rangle}}

\newcommand{\CS}{{\mathcal{S}}}

\newcommand{\CJ}{{\mathcal{J}}}
\newcommand{\wedgeq}{{\wedge\kern-5pt\cdot}}

\newcommand{\tens}{\otimes}

\newcommand{\id}{{\rm id}}

\newcommand{\extd}{{\rm d}}
\newcommand{\del}{{\partial}}
\newcommand{\eps}{\epsilon}

\newcommand{\la}{{\triangleright}}
\newcommand{\ra}{{\triangleleft}}

\newcommand{\dirac}{{ \slashed{D} }}

\newcommand{\doublenabla}{%
  \nabla\mkern-12mu\nabla}   
  
\begin{document}

\title{Dirac operator associated to a quantum metric}


\author{S. Majid}
\address{Queen Mary, University of London\\
School of Mathematics, Mile End Rd, London E1 4NS, UK}

\email{s.majid@qmul.ac.uk}


\begin{abstract}  We construct a canonical geometrically realised Connes spectral triple or `Dirac operator' $\dirac$ from the data of a quantum metric $g\in \Omega^1\tens_A\Omega^1$ and quantum Levi-Civita bimodule connection, at the pre-Hilbert space level. Here $A$ is a possibly noncommutative coordinate algebra, $\Omega^1$ a bimodule of 1-forms and the spinor bundle is $\CS=A\oplus\Omega^1$. When applied to graphs or lattices, we essentially recover a previously proposed Dirac operator but now as a geometrically realised spectral triple. We also apply the construction to the fuzzy sphere and to $2\times 2$ matrices with their standard quantum Riemannian geometries. We also propose how $\dirac$ can be extended to an external bundle with bimodule connection. \end{abstract}
\maketitle 

\section{Introduction}

One of the important ideas for noncommutative geometry is an approach coming out of cyclic cohomology and operator algebras\cite{Con} leading up to the notion of a `Connes spectral triple' in the role of Dirac operator\cite{Con0,Con:spe}. Essentially, this is: 

(1) a Hilbert space $\CS,\<\ ,\ \>$ on which a $*$-algebra $A$ acts faithfully with the action of $a^*$ adjoint to the action of $a$, and on which operators $\dirac,\CJ$ and in the `even' case $\gamma$ act, with $\dirac$  (in our conventions) antihermitian, $\CJ$  an antilinear isometry and $\gamma$ hermitian. 

(2) $[a,\CJ b \CJ^{-1}]=[[\dirac,a], \CJ b\CJ^{-1}]=[\gamma, a]=0$ for all $a,b\in A$.

(3) $\CJ^2=\eps\id$, $\CJ\dirac =\eps'\dirac \CJ$, $\CJ\gamma=\eps''\gamma\CJ$, $\gamma^2=\id$, $\dirac\gamma=-\gamma\dirac$.

for signs $\eps,\eps',\eps''$ which by analogy with the classical case are usually taken to fall into a period 8 table according to a `dimension' $n$. Here $\eps'$ is opposite to Connes' conventions since we take $\dirac$ antihermitian to match our algebraic formulation. In this paper we will {\em only} construct spectral triples with $\eps=\eps'=\eps''=1$, which do not then fully fit into this table: one could consider them $n=1$ ignoring $\gamma$ from this point of view. If we solve for just (2)-(3) without a Hilbert space structure on $\CS$ then we say that we have a spectral triple at the {\em local tensorial level}. After achieving the Hilbert space level, one usually demands additionally that $\dirac$ has compact resolvent and bounded commutators. 

In classical geometry, the local tensorial level would be the construction of $\dirac$, charge conjugation $\CJ$ and a chirality operator $\gamma$ at the level of smooth sections of the spinor bundle $\CS$ before completing to $L^2$ sections. In this paper,  we continue the programme in \cite{BegMa:spe,LirMa2,LirMa3}\cite[Chapter~8.5]{BegMa} to `quantum-geometrically realise' Connes spectral triples starting first at the local tensorial level and in such a way that $\dirac=\la\circ\nabla_S$, where $\nabla_S:\CS\to \Omega^1\tens_A\CS$ is a bimodule connection on an $A$-bimodule $\CS$ and $\la:\Omega^1\tens_A\CS\to \CS$ is a bimodule map (the `Clifford action') in the role of gamma matrices. We use here a more `layer by layer' approach to noncommutative geometry in which the starting point is an $A$-bimodule $\Omega^1$ of `differential forms' equipped with a derivation $\extd: A\to \Omega^1$. In recent years, this was extended to a systematic theory of  `quantum Riemannian geometry' (QRG) with a quantum metric $g\in \Omega^1\tens_A\Omega^1$, a `quantum Levi-Civita connection' (QLC) $\nabla:\Omega^1\to \Omega^1\tens_A\Omega^1$, etc. We refer to \cite{BegMa} and references therein. This framework has been applied to baby models of quantum gravity and quantum geodesics, and is also beginning to be applied to particle physics, see \cite{ArgMa4} for a recent review. The geometric realisation programme is about the intersection of these two approaches to noncommutative geometry. 

In this context, we provide and study a specific construction that associates a geometrically realised spectral triple on the spinor bundle $\CS=A\oplus\Omega^1$ to every QRG under some mild conditions. If not all of these conditions hold then we can still obtain a natural $\dirac$ at the local tensorial level obeying all or most of (2)-(3). Our general construction appears in Section~\ref{secdir}, see Theorem~\ref{thmspe}. We then compute how this works in three examples, the first of which in Section~\ref{secgra} produces results broadly similar to the Dirac operators studied on graphs studied several times since \cite{Dav}, most recently in \cite{Bia}.  In all three examples the calculus is inner and in this case, Corollary~\ref{baredir}, there is a canonical but strictly quantum `bare' choice of $\nabla$ where $\sigma=\id$ so that the only input is the quantum metric. Remarkably, in the examples that we compute, the Dirac operator is the same as obtained by using the QLC. And because of the general context, we can also take the classical limit of our construction, where it turns out to be the projection to degrees 0,1 of the spectral triple in \cite[Thm.~11.3]{Con:spe} on any closed orientable Riemannian manifold $M$ with $\CS=\Omega(M)_\C$ the complexified exterior algebra and $\dirac=\extd+\delta$ for the Hodge co-differential $\delta$. Connes showed in this work that spectral triples in the commutative case, subject to some other mild assumptions, are either of this type or given by a usual Dirac operator on spinors. On the other hand, there is no straightforward Hodge-de Rham theory in the general noncommutative case (although these are reasonable examples in special cases \cite{Ma:hod, Ma:rec}), whereas our degrees 0,1 version will apply at the local tensorial level for any QRG with a certain $\sigma$-symmetry condition on the quantum metric.  

It is worth noting that $\CS=A\oplus\Omega^1$ is also the first jet bundle $\CJ^1$ in \cite{MaSim} but while the jet bundle extends to higher order in a quantum symmetric manner, $\CS$ is more like the leading part of the exterior algebra which classically would be  antisymmetric. The same applies when we extend the construction to include an external bundle $E$ with connection $\nabla_E$, which we briefly consider in Section~\ref{secE}. This extended construction in the graph case seems somewhat different from the lattice Dirac operator coupled to an external gauge field in \cite{Bia1} but has elements in common. The paper concludes with some remarks in Section~\ref{secrem} about directions for further work.

We work over $\C$ but all of the constructions at the local tensorial level work over any field, other than those involving $\CJ$. For axioms (1), we actually work at a pre-Hilbert space level without consideration of any analytic issues needed for extension to the $L^2$ completion, compact resolvent and bounded commutators (what could more precisely be called a pre-spectral triple). This suffices for our purposes, where examples will be finite-dimensional or essentially so. 
  
  \subsection*{Acknowledgements} I would like to thank G. Bianconi for explaining to me her papers \cite{Bia,Bia1} and references therein for the graph case, which was a motivation behind the present work. 
 
\section{QRG's and the construction of a spectral triple}

We start with some preliminaries on the QRG formalism and on how this can be used to construct geometric spectral triples\cite{BegMa, BegMa:spe}. Then Section~\ref{secdir} contains our general construction within this formalism  given a bimodule connection $\nabla,\sigma$, a $\sigma$-symmetric quantum metric and a suitable integration state. 

\subsection{Elements of QRG}\label{secele} We work with $A$ a unital algebra, typically a $*$-algebra over $\C$, in the role of `coordinate algebra'. Differentials are formally introduced as a bimodule $\Omega^1$ of 1-forms equipped with a map $\extd:A\to \Omega^1$ obeying the Leibniz rule $\extd(ab)=(\extd a)b+a\extd b$. We assume this extends to an exterior algebra $(\Omega,\extd)$ with $\extd^2=0$ and $\extd$ obeying the graded-Leibniz rule, and where $A,\extd A$ generate $\Omega$ (this is more restrictive than a differential graded algebra). 

A quantum metric is $g\in \Omega^1\tens_A\Omega^1$ and a bimodule map inverse $(\ ,\ ):\Omega^1\tens_A\Omega^1\to A$ in an obvious sense, together with some form of quantum symmetry condition.  Following \cite{BegMa}, we say $g$ is a generalised quantum metric if no form of symmetry is imposed. A (left) bimodule connection\cite{DVMic,Mou} on $\Omega^1$ is $\nabla:\Omega^1\to \Omega^1\tens_A\Omega^1$ obeying 
\begin{equation}\label{nabla} \nabla(a \omega)=a \nabla\omega+ \extd a\tens\omega,\quad \nabla(\omega a)=(\nabla\omega) a+\sigma(\omega\tens\extd a)\end{equation}
for all $a\in A,\omega\in \Omega^1$, for some `generalised braiding' bimodule map $\sigma:\Omega^1\tens_A\Omega^1\to \Omega^1\tens_A\Omega^1$. The latter, if it exists, is uniquely determined and not additional data. A connection is  {\em metric compatible} if the tensor product connection
\[ \nabla_{\Omega^1\tens
\Omega^1}:=\nabla\tens \id+ (\sigma\tens\id)(\id\tens\nabla)\]
vanishes on $g$. Provided $\sigma$ is invertible (which we assume), metric compatibility is equivalent to covariance of $(\ , \ )$ in the sense
\begin{equation}\label{invmetcomp} \extd(\ ,\ )=(\id\tens (\ ,\ ))\nabla_{\Omega^1\tens\Omega^1}    \end{equation}
as shown in \cite[Lem.~8.4]{BegMa}. It is this form which we will use (we do not actually need $g$ itself).  In this paper, the natural symmetry condition for a quantum metric we will be led to impose is `$\sigma$-symmetric'
\begin{equation}\label{sig-sym} (\ ,\ )\sigma=(\ ,\ )\end{equation}
 in terms  of $(\ ,\ )$, with respect to a bimodule connection $\nabla,\sigma$. This is less self-contained than the standard notion $\wedge(g)=0$ of quantum symmetry in \cite{BegMa}.  
 
 When $\nabla$ is both metric compatible and torsion free in the sense that the torsion $T_\nabla:=\wedge\nabla-\extd$ vanishes, we say that we have a {\em quantum Levi-Civita connection} (QLC). In the torsion free case, one can show that $\wedge\circ(\id+\sigma)=0$. 

Finally, working over $\C$, we need  $(\Omega,\extd)$ to be a $*$-calculus, where $A$ is a $*$-algebra and $\extd$ is extended so as to commute with $*$. We also usually require $g$ `real' in the sense $g^\dagger=g$, where $\dagger={\rm flip}(*\tens *)$ is well defined on $\Omega^1\tens_A\Omega^1$, or equivalently 
\begin{equation}\label{realmet} *\circ(\ ,\ )=(\ ,\ )\circ\dagger\end{equation}
in terms of the inverse metric. Moreover, we usually require $\nabla$ $*$-preserving in the sense 
\begin{equation}\label{nabla*}\nabla\circ *=\sigma\circ\dagger\circ\nabla.\end{equation}
This in turn can be shown\cite[Lemma~8.7]{BegMa} to imply that $(\sigma\circ\dagger)^2=\id$. 

We will particularly need in this context a natural co-differential or divergence 
\[ \delta={\rm div}_\nabla:\Omega^1\to A,\quad \delta :=(\ ,\ )\nabla\]
which necessarily obeys 
\begin{equation}\label{delta} \delta(a\omega)=(\extd a,\omega)+a\delta\omega,\quad \delta(\omega a)=(\ ,\ )\sigma(\omega\tens \extd a)+(\delta\omega)a\end{equation}
 for all $a\in A, \omega\in \Omega^1$. 
 
 \begin{lemma}\label{delta*} In the $*$-calculus case, if the metric and $\nabla$ are obey their respective reality and $*$-preserving conditions and the metric is $\sigma$-symmetric then $\delta$ commutes with $*$. 
\end{lemma}
\proof Under our assumptions, $(\delta\omega)^*=(\ ,\ )\dagger \nabla\omega=(\ ,\ )\sigma^{-1}\nabla(\omega^*)=\delta(\omega^*)$ for all $\omega\in\Omega^1$.
\endproof

 This also implies that the Laplacian $\Delta=(\ ,\ )\nabla\extd=\delta\extd$ on $A$ commutes with $*$ as in \cite{BegMa}.  Following \cite{BegMa:cur}, we say that  a positive-linear functional $\int: A\to \C$ is {\em divergence-compatible} if 
\begin{equation}\label{div-compat} \int \circ\, \delta=0\end{equation}
as would be true in the classical case for a closed manifold. In this case (\ref{delta}) integrates to 0 for each left hand side, giving formulae for integration by parts.

\subsection{Spinor bundle and Dirac operator} \label{secspe}

The formulation of a bimodule connection above applies equally well to any $A$-bimodule  $\CS$ in the role of sections of a `vector bundle', now with $\nabla_S:\CS\to \Omega^1\tens_A\CS$ and $\sigma_S:\CS\tens_A\Omega^1\to \Omega^1\tens_A\CS$, and (\ref{nabla}) now with $\omega\in \Omega^1$ replaced by $
\phi\in \CS$. We also need a {\em Clifford bimodule map} $\la:\Omega^1\tens_A\CS\to \CS$ which generalises the role of the gamma-matrices when constructing Dirac operators\cite{BegMa:spe}. In this case, we define $\dirac=\la\circ\nabla_S: \CS\to \CS$ as a minimal `geometric Dirac operator' in the formalism. 

Not required for a spectral triple but motivated by the classical case, we optionally would like the geometric condition $\doublenabla(\la)=0$ i.e., 
\begin{equation}\label{covla}  (\id\tens\la)\circ(\nabla\tens \id+ (\sigma\tens\id)(\id\tens\nabla_S))=\nabla_S\circ\la,\end{equation}
which says that $\la$ is covariantly constant in the sense that it intertwines the tensor product connection on $\Omega^1\tens_A\CS$ and the connection on $\CS$. \cite{BegMa:spe} also tentatively proposed a `Clifford relation' that  $\la$ extends to $\Omega^2\tens_A\CS\to \CS$ in such a way that 
\begin{equation}\label{cliffla} (\omega\wedge\eta)\la \phi:=\omega\la(\eta\la s)+(\omega,\eta)\phi\end{equation}
for all $\omega,\eta\in \Omega^1$ and $\phi\in \CS$, is well defined. One can also consider a generalisation of this to allow an automorphism of $\CS$ and a scaling of the inner product. The motivation is that this reduces to the Clifford algebra relations in the classical case. This  condition, however, is less natural and does not fully apply in many examples, including in the constructions below. It is also not needed for a spectral triple. 

For a Connes spectral triple, we instead need\cite{BegMa:spe}\cite[Chapter~8.5]{BegMa} that $\CS$ has an antilinear skew-bimodule `charge conjugation' involution $\CJ:\CS\to \CS$  in the sense
\begin{equation}\label{J} \CJ(a\phi)=\CJ(\phi)a^*,\quad \CJ(\phi a)=a^*\CJ(\phi),\quad  \CJ^2=\eps\id\end{equation}
for all $a\in A$, $\phi\in \CS$ and some $\eps=\pm1$. This is required to obey
\begin{equation}\label{Jeqns} \nabla_S\CJ= \sigma_S\circ{\rm flip}(*\tens \CJ)\circ\nabla_S,\quad \CJ(\omega\la\phi)=\eps'\la\circ\sigma_S(\CJ\phi\tens\omega^*). \end{equation}
The first equation generalises the notion of $*$-preserving for $\nabla$ to $\nabla_S$. Also, for an even spectral triple, we need a bimodule map involution  $\gamma:\CS\to \CS$ with $\gamma^2=\id$ required to further obey
\begin{equation}\label{gam}   (\id\tens\gamma)\circ\nabla_S=\nabla_S\circ\gamma, \quad \gamma\circ\la=-\la\circ(\id\tens\gamma), \quad \CJ\circ\gamma=\eps''\gamma\circ\CJ.\end{equation}
It is shown in \cite{BegMa:spe} that in this way we naturally realise the local tensorial part of Connes axioms of a spectral triple.  

Finally, we need a compatible sesquilinear inner product $\<\ ,\ \>$ on $\CS$ with respect to which we complete the latter into a Hilbert space as stated in the introduction, namely that  $\dirac$ is antihermitian, $\CJ$ an antilinear isometry in the sense
\[ \<\CJ\phi,\CJ\psi\>=\<\psi,\phi\>\]
and (if it exists) $\gamma$ hermitian. Here, $\imath\dirac$ is the hermitian Dirac operator in Connes axioms with the result that our $\eps'$ has the reversed sign compared to Connes conventions. Our approach will be to  build  $\<\ ,\ \>$ from a divergence-compatible positive linear trace $\int:A\to \C$. 

\subsection{Construction of $\dirac$  from a QRG}
\label{secdir}

We are now ready to give our construction. Given a differential algebra $(A,\Omega^1,\extd)$, we set
\begin{equation}\label{CS} \CS=A\oplus \Omega^1\end{equation}
 with its direct sum bimodule structure. We next assume a generalised quantum metric $g$ with inverse $(\ , \ )$ and $\nabla$ a bimodule connection on $\Omega^1$, and define the bimodule connection
\begin{equation}\label{nablaS} \nabla_S (a+\omega)=\extd a+ \nabla\omega\end{equation}
as an element of $\Omega^1\tens_A\CS=\Omega^1\oplus \Omega^1\tens_A\Omega^1$ after the canonical identification $\Omega^1\tens_A A=\Omega^1$. Here, as elsewhere,  $a\in A$ and $\omega,\eta\in \Omega^1$. The associated generalised braiding is 
\begin{equation}\label{sigS} \sigma_S((a+\omega)\tens\eta))=a\eta+ \sigma(\omega\tens\eta)\end{equation}
where also, on the input side, $\CS\tens_A\Omega^1=\Omega^1\oplus\Omega^1\tens_A\Omega^1$ by the canonical identification $A\tens_A\Omega^1=\Omega^1$. 
Next, for the Clifford action, we define
\begin{equation}\label{la} \omega\la (a+\eta)=\omega a+(\omega,\eta)\end{equation}
so that 
\begin{equation}\label{dir}\dirac(a+\omega)=\la\nabla_S(a+\omega)=\extd a+(\ ,\ )\nabla\omega\end{equation} or in other words $\dirac=\extd+\delta$ 
as desired. From the Leibniz rules, this necessarily obeys 
\begin{equation}\label{comdir}[\dirac,a]=(\extd a)\la\end{equation}
for all $a\in A$. Similarly, from (\ref{dir}) it follows that
\begin{equation}\label{dirsq}\dirac^2(a+\omega)=(\ ,\ )\nabla\extd a+ (\id\tens (\ ,\ ))\nabla_{\Omega^1\tens\Omega^1}
\nabla \omega=\Delta a+\Delta\omega\end{equation}
as the natural Laplacians on functions and 1-forms  in \cite{Ma:gra}\cite[Lem.~8.6]{BegMa}. 

\begin{lemma} $\doublenabla(\la)=0$ in (\ref{covla}) holds iff $\nabla$ is metric compatible.
\end{lemma}
\proof  For covariance of $\omega\tens a\mapsto \omega a$ we need
\[ \nabla(\omega a)=(\id\tens\la)(\nabla\omega\tens a+\sigma(\omega\tens \extd a)=(\nabla \omega)a+\sigma(\omega\tens\extd a)\]
so this recovers the condition for $\nabla$ to be a bimodule connection. For covariance of $\omega\tens\eta\mapsto (\omega,\eta)$ we need
\[ \extd(\omega,\eta)=(\id\tens\la)\left(\nabla\omega\tens \eta+(\sigma\tens\id)(
\omega\tens\nabla\eta)\right)=(\id\tens (\ ,\eta))\nabla\omega+(\id\tens (\ ,\ ))(\sigma\tens\id)(\omega\tens\nabla\eta)\]
which is equivalent to metric compatibility when written in terms of $(\ ,\ )$ (see \cite[Chap. 8]{BegMa}). 
 \endproof

Note that we will not actually need $T_\nabla=0$, so $\nabla$ does not need to be a full QLC, although this would be the natural choice if it exists. On the other hand, the Clifford relations (\ref{cliffla}) become
\[ \omega\la (\eta\la a)- (\omega,\eta)a=0=:(\omega\wedge\eta)\la a\]
\[ \omega\la(\eta\la\zeta)-(\omega,\eta)\zeta=\omega(\eta,\zeta)-(\omega,\eta)\zeta=:(\omega\wedge\eta)\la\zeta\]
which seems too strong in that the 2nd equation would not be well defined classically as the wedge product is antisymmetric in $\omega,\eta$. Hence, we will not try to impose this (optional) condition. Next, we define
\begin{equation}\label{Jex} \CJ(a+\omega)=a^*+\omega^*\end{equation}
which automatically obeys (\ref{J}) with $\eps=1$. 
\begin{lemma}\label{lemJ} Suppose that the generalised quantum metric obeys the reality (\ref{realmet}) and $\nabla$ is $*$-preserving. Then $\CJ$ obeys  (\ref{Jeqns}) iff the metric is  $\sigma$-symmetric and  $\eps'=1$.
\end{lemma}
\proof  For the first of (\ref{Jeqns}), we have
\begin{align*} \sigma_S{\rm flip}(*\tens\CJ)(\extd a+\nabla\omega)&=\sigma_S(1\tens\extd a^*+\dagger\circ\nabla\omega)=\sigma_S(1\tens\extd a^*+\sigma^{-1}\circ\nabla(\omega^*))\\
&=\extd a^*+\nabla(\omega^*)=\nabla_S(a^*+\omega^*)=\nabla_S\CJ(a+\omega)\end{align*}
so that this holds automatically given that $\extd$ commutes with $*$ and $\nabla$ is $*$-preserving. For the second of (\ref{Jeqns}) we have
\begin{align*}\eps'\la\circ\sigma_S(\CJ(a+\omega)\tens\eta^*)&=
\eps'\la\circ\sigma_S((a^*+\omega^*)\tens\eta^*)=\eps'\la(a^*\eta^*\tens 1+\sigma(\omega^*\tens\eta^*))\\
&=\eps'(a^*\eta^*+(\ ,\ )\sigma(\omega^*\tens\eta^*))
\end{align*}
and we want this to equal 
\[\CJ(\eta\la(a+\omega))=\CJ(\eta a+(\eta,\omega))=(\eta,\omega)^*+a^*\eta^*=a^*\eta^*+(\omega^*,\eta^*)\]
using the reality of the metric in the $(\ ,\ )$ form as in \cite[Chap 1]{BegMa}. Comparing, we need the conditions stated. 
\endproof

In fact, what we actually used in the proof (and this is all we need for a spectral triple) is neither usual metric reality nor $\sigma$-symmetry but  their combination
\begin{equation}\label{Jmet} (\omega,\eta)^*=(\ ,\ )\circ\sigma(\eta^*\tens\omega^*),\end{equation}
 which is  the minimum condition needed on the generalised quantum metric for the present construction to work. Finally, we set
\begin{equation}\label{gamS} \gamma(a+\omega)=a-\omega,\end{equation}
which is clearly an involution and a bimodule map, and clearly obeys the last part of (\ref{gam}) with $\eps''=1$.

\begin{lemma}\label{lemgam} $\gamma$ obeys the rest of (\ref{gam}) also.
\end{lemma}
\proof Here the first part of (\ref{gam}) is
\[ (\id\tens\gamma)\nabla_S(a+\omega)=(\id\tens\gamma)(\extd a\tens 1+\nabla\omega)=\extd a-\nabla\omega=\nabla_S(a-\omega)=\nabla_S\gamma(a+\omega)\]
and the 2nd part is
\[ \gamma(\eta\la(a+\omega))=\gamma(\eta a+(\eta,\omega))=-\eta a-(\eta,-\omega)=-\eta\la \gamma(a+\omega)\]
as required. 
 \endproof

If the conditions in these lemmas hold then we have all the local tensorial axioms of a spectral triple. It remains to define the inner product and we do this by
\begin{equation}\label{innS} \<a+\omega, b+\eta\>=\int (a^* b+ (\omega^*,\eta))\end{equation}
which is manifestly sesquilinear as needed provided $\int$ is a *-preserving linear functional and the metric reality (\ref{realmet}) holds at least under the integration. We assume $\int$ is {\em extended-positive definite} in the sense
\begin{equation}\label{extposdef} \int(a^*a)\ge 0,\quad \int(\omega^*,\omega)\ge 0\end{equation}
for all $a\in A,\omega\in\Omega^1$ with equality iff $a=0$ or $\omega=0$ respectively. These assumptions make $\CS$ into a pre-Hilbert space which can then be completed  to a Hilbert space. 

\begin{theorem}\label{thmspe} If the quantum metric obeys the reality condition (\ref{realmet}) and $\int$ is extended-positive definite and divergence-compatible then $\dirac$ above is antihermitian and $\gamma$ is hermitian. Moreover, $\CJ$ is an antilinear isometry iff $\int$ is an {\em extended trace} in the sense
\[ \int(ab)=\int(ba), \quad \int (\omega,\eta)=\int (\eta,\omega)\]
for all $a,b\in A$ and $\omega,\eta\in \Omega^1$. 
In this case, we have a (pre)spectral triple with $\eps=\eps'=\eps''=1$.
\end{theorem}
\proof Here, using integration by parts
\begin{align*} \<a+\omega,\dirac(b+\eta)\>&=\<a+\omega, \extd b+\delta\eta\>=
\int (a^*\delta\eta+(
\omega^*,\extd b))=\int(-(\extd a^*,\eta)+(
\omega^*,\extd b))\\
\<-\dirac(a+\omega),b+\eta\>&=\<-\extd a-\delta\omega,b+\eta\>=-\int((\delta\omega)^*b+(\extd a^*,\eta))\\
&=-\int((b^*\delta\omega)^*+(\extd a^*,\eta))=\int((\extd b^*,\omega)^*-(\extd a^*,\eta))
\end{align*}
where for the last equality we used that $\int$ is $*$-preserving to allow us to kill the $*$ of a total divergence. We then need (\ref{realmet}) to hold at least under the integral to obtain the same result as the first calculation. Hence $\dirac$ is antihermitian. Next, 
\[ \<\gamma(a+\omega),\gamma(b+\eta)\>=\<a-\omega,b-\eta\>=\<a+\omega,b+\eta\>\]
is immediate from the formula for the inner product. Finally,
\begin{align*}\<\CJ(a+\omega),\CJ(b+\eta)\>&=\<a^*+\omega^*,b^*+\eta^*\>=\int(a b^*+(\omega,\eta^*))\\
\<b+\eta,a+\omega\>&=\int(b^*a+(\eta^*,\omega))
\end{align*}
Equality needs the `extended trace' conditions stated. \endproof

 From the proof, we only really need the reality of the metric under the $\int$ rather than (\ref{realmet}). Remembering the sign change of $\eps'$ when passing to the hermitian $\imath\dirac$, the spectral triple here is technically $n=1$ ignoring $\gamma$ in the Connes period 8 classification table for the pattern of signs. 

Also recall that in the local part of the construction we did not actually need $\nabla$ to be metric compatible as this was only used for the optional (\ref{covla}). However, if there is a natural $\nabla$ which is metric compatible (such as a QLC) then it makes sense to use this and consider the above $\nabla_S$ as the `base' coming from the QRG. Then any other bimodule connection on $\CS$ with the same $\sigma_S$ differs by a bimodule map of the form
\begin{equation}\label{alphaS} \alpha_S:\CS\to \Omega^1\tens_A\CS,\quad 
  \alpha_S(a+\omega)=a\alpha_0+\alpha(\omega)\end{equation}
for some central 1-form $\alpha_0$ and some bimodule map $\alpha:\Omega^1\to \Omega^1\tens_A\Omega^1$. We can then  use the bimodule connection $\nabla_{\alpha_S}:=\nabla_S+\alpha_S$ in place of $\nabla_S$ giving a modified Dirac operator 
\[\dirac_{\alpha_S}=\dirac+ \la\circ\alpha_S,\quad  (\la\circ\alpha_S)(a+\omega)=a\alpha_0+(\ ,\ )\alpha(\omega).\]
\begin{proposition}\label{a0alpha} The equations (\ref{Jeqns})-(\ref{gam}) continue to hold for $\nabla_{\alpha_S}$ provided
\[ \alpha_0^*=\alpha_0,\quad \alpha\circ *=\sigma\circ\dagger\circ\alpha \]
and in this case we still have a (pre)spectral triple  with $\dirac_{\alpha_S}$ provided
\[ (\ ,\ )\alpha+(\alpha_0,\ )=0.\]
\end{proposition}
\proof We have to recheck the proofs involving $\nabla_S$ when $\alpha_S$ is added. For (the first part of) (\ref{Jeqns}) we need 
\[ \sigma_S{\rm flip}(*\tens\CJ)(a\alpha_0\tens 1+\alpha(\omega))=\sigma\circ\dagger\circ\alpha(\omega)+\alpha_0^* a^*\]
to equal $\alpha_S\CJ(a+\omega)=a^*\alpha_0+\alpha(\omega^*)$, which gives the conditions stated given that $\alpha_0$ is central. For the (first part of) (\ref{gam}) we have 
\[ (\id\tens\gamma)(a\alpha_0\tens 1+\alpha(\omega))=a\alpha_0-\alpha(\omega)=\alpha_S(a-\omega)=\alpha_S\gamma,\]
as required, so this is automatic. For the extension of Theorem~\ref{thmspe} it remains to check that the extended $\dirac_{\alpha_S}$ is anti-hermitian. Here
\begin{align*} \<a+\omega,\la\circ\alpha_S(b+\eta)\>&=\<a+\omega, b\alpha_0+(\ , )\alpha(\eta)\>=\int(a^*(\ ,\ )\alpha(\eta)+ (\omega^*,\alpha_0)b)\\
\<-\la\circ\alpha_S(a+\omega),b+\eta\>&=\<-a\alpha_0-(\ ,\ )\alpha(\omega),b+\eta\> =\int(-a^*(\alpha_0^*,\eta)-((\ ,\ )\alpha(\omega))^*b)\end{align*}
using that $\alpha_0$ is central. Equality is ensured, given that $\alpha_0^*=\alpha_0$ from the earlier part and given that the metric obeys the reality condition, provided we assume the condition on $(\ ,\ )\alpha$ stated (and this is required if $\int$ is suitably nondegenerate). \endproof

 The last condition here says that $\alpha_0$ is determined from the bimodule map $\alpha$, so the latter is the additional freedom for the Dirac operator constructed on the same $\CS$ and for given $\sigma_S$  (other connections and hence Dirac operators could be possible e.g. using a connection for $\nabla$ with a different $\sigma$). Although this represents additional freedom in the construction of $\dirac$ and has the nature of a gauge field, it is not the most general geometric way to introduce an external gauge symmetry, which should be done by tensoring with an external bundle with connection as we do in Section~\ref{secE}. 
 
\section{Graph, matrix and fuzzy Dirac operators}\label{secex}

Here we analyse how to solve the conditions needed for our construction in some important cases, the first of which is quite general but has no classical analogue. 

\subsection{Bare Dirac operator for an inner calculus}

A calculus $(\Omega^1,\extd)$  is inner if there is $\theta\in \Omega^1$ such that $\extd a=[\theta,a]$. This is not possible in classical geometry but quite typical in the noncommutative case. In this case, bimodule connections have the form\cite{Ma:gra}
\[ \nabla\omega=\theta\tens\omega-\sigma(\omega\tens\theta)+\alpha(\omega)\]
for arbitrary bimodule maps $\sigma:\Omega^1\tens_A
\Omega^1\to \Omega^1\tens_A
\Omega^1$ and $\alpha:\Omega^1\to \Omega^1\tens_A\Omega^1$. If there is no such last term, we say that $\nabla$ is inner. The choice $\sigma=s\, \id$ is a bimodule map for any constant $s$, giving  a 1-parameter canonical inner connection on any inner calculus. Here 
\[ \nabla\omega=\theta\tens\omega-s\omega\tens\theta\]
so that the corresponding Laplacian is
\[ \Delta a=(\ ,\ )\nabla\extd a=(\theta,\extd a)-s(\extd a,\theta)=(\theta,\theta)a+s a(\theta,\theta)-(1+s)(\theta,a\theta).\]
The choice $s=0$ is (up to a factor of 2) the Laplacian $\Delta_\theta$ in \cite[Chap 1.5]{BegMa}. The choice $s=1$ by contrast is more suited to our needs and trivially obeys the $\sigma$-symmetric condition in Lemma~\ref{lemJ} for any generalised quantum metric. The bare $\nabla$ is also $*$-preserving if 
\[ \theta^*=-\theta\]
 which we assume. This bare $\nabla$ therefore generally leads to a Dirac operator and spectral triple at the local or tensor level. 

\begin{corollary}\label{baredir} If $\int$ is an extended trace as in Theorem~\ref{thmspe} then it is automatically divergence-compatible. In this case we have {\em bare dirac operator}
\[ \dirac (a+\omega)=\extd a+  (\theta,\omega)-(\omega,\theta)\]
forming a geometrically realised (pre)spectral triple with $\eps=\eps'=\eps''=1$ if $\int$ is extended-positive definite on $\CS$. 
\end{corollary} 
\proof Here $\delta\omega=(\ ,\ )\nabla\omega=(\theta,\omega)-(\omega,\theta)$ 
so if the integral obeys the integral symmetry condition then it is automatically divergence-compatible. We then apply Theorem~\ref{thmspe}. To have an pre-Hilbert space we still need $\int$ to be extended-positive definite as in (\ref{extposdef}). \endproof

In addition, we can add non-inner terms $\alpha$ to the bare case above. On the other hand, we will generally not have the optional $\doublenabla(\la)=0$ covariance condition as the `bare' $\nabla$ itself is not generally metric compatible, being usually very far from a quantum Levi-Civita connection. 

\subsection{Graph case}\label{secgra} As in \cite{Ma:gra}
\cite[Chap 1.5]{BegMa}, we can take $A=\C(X)$ the complex functions on the vertex set $X$ of a bidirected graph (i.e., every edge has arrows in both directions) and $\Omega^1$ with vector space basis $\omega_{x\to y}$ labelled by the arrows. Then spinor space in our construction is $\CS=\C(X)\oplus\Omega^1$, which in the finite graph case is $|X|+2|E|$-dimensional over the field. We focus on the finite case but the underlying constructions work more generally with suitable care.  

The bimodule structure of $\Omega^1$ and $\extd$ are well-known, 
\[ a.\omega_{x\to y}=a(x)\omega_{x\to y},\quad \omega_{x\to y}.a=a(y)\omega_{x\to y},\quad \extd a=\sum_{x\to y}(a(y)-a(x))\omega_{x\to y},\quad \omega_{x\to y}^*=-\omega_{y\to x}.\]
A (generalised) quantum metric has the form
\[ (\omega_{x\to y},\omega_{y'\to x'})=\lambda_{x\to y}\delta_{x,x'}\delta_{y,y'}\delta_x \]
for real coefficients $\lambda_{x\to y}$ associated to the arrows. It is typical in this context (but not necessary) to assume that the metric is {\em edge-symmetric} in the sense $\lambda_{x\to y}=\lambda_{y\to x}$ so that we depend only on the edge, but we do not assume this unless stated. Elements of $\Omega^1\tens_A\Omega^1$ have vector space basis $\{\omega_{x\to y\to z}\}$ labelled by the 2-steps. We let
\[ \Omega^2_{x,y}={\rm span}_\C\{ \omega_{x\to w\to y}\}\]
be the 2-steps with fixed endpoints with basis  labeled by the possible intermediate vertices. 

Among the bimodule connections $\nabla$,  the simplest are of  `inner form'
\[ \nabla \omega=\theta\tens\omega -\sigma(\omega\tens\theta),\quad\theta=\sum_{x\to y}\omega_{x\to y}\]
for any bimodule map $\sigma:\Omega^1\tens_A\Omega^1\to \Omega^1\tens_A\Omega^1$.  If there are no triangles in the graph then there is no possibility for a bimodule map $\alpha$ and hence all connections are inner. Moreover,  $\sigma$ amounts to a free choice of  linear maps on each $\Omega^2_{x,y}$ space, i.e. to matrices in $\sigma_{x,y}\in M_{\dim(\Omega^2_{x,y})}(\C)$ for all pairs $x,y$ that can be endpoints of 2-steps. One usually restricts to invertible $\sigma$, which means each of these matrices should be invertible. In these terms, an inner connection looks like
\[ \nabla(\omega_{x\to y})=\sum_{w:w\to x}\omega_{w\to x\to y}-\sum_{w:y\to w}\sigma(\omega_{x\to y\to w}) =\sum_{w:w\to x}\omega_{w\to x\to y}-\sum_{w:y\to w}\sum_{z: x\to z\to w}\sigma_{x,w}{}^z{}_y \omega_{x\to z\to w}\]
and is $*$-preserving if $\dagger\sigma$ squares to the identity, which amounts to 
\begin{equation}\label{siggraph*} \sigma_{y,x}=\overline{\sigma_{x,y}^{-1}}\end{equation}
as matrices. Here $\Omega^2_{x,y}$ is identified with $\Omega^2_{y,x}$ by arrow reversal of the basis elements. The divergence $\delta=(\ ,\ )\nabla$ for the chosen bimodule connection  is then
\begin{equation}\label{gradiv} \delta(\omega_{x\to y})=\lambda_{y\to x}\delta_y-\left(\sum_{z:x\to z}\lambda_{x\to z}\sigma_{x,x}{}^z{}_y\right)\delta_x  \end{equation}
which is sensitive only to the $\sigma_{x,x}$ matrices. The somewhat general quantum geometric Dirac operator by our construction (but not necessarily yet a spectral triple)  is therefore
\begin{equation}\label{gradir}\dirac(\delta_x)=\sum_{y:y\to x}\omega_{y\to x}- \sum_{y: x\to y}\omega_{x\to y},\quad \dirac(\omega_{x\to y})= \lambda_{y\to x}\delta_y-\left(\sum_{z:x\to z}\lambda_{x\to z}\sigma_{x,x}{}^z{}_y\right)\delta_x\end{equation}

Next, towards a spectral triple, we need associated maps $\CJ$ and $\gamma$, which in our construction are
\begin{equation}\label{graJ} \CJ(\delta_x)=\delta_x,\quad \CJ(\omega_{x\to y})=-\omega_{y\to x},\quad \gamma(\delta_x)=\delta_x,\quad \gamma(\omega_{x\to y})=-\omega_{x\to y}\end{equation}
extended antilinearly and linearly respectively. The $\sigma$-symmetric condition needed in Lemma~\ref{lemJ} is
\begin{equation}\label{grasym} \lambda_{x\to y}=\sum_{z:x\to z}\lambda_{x\to z}\sigma_{x,x}{}^z{}_y\end{equation}
for all $x\to y$, and completes the requirements for a geometrically realised spectral triple at the local tensorial level. The Dirac operator is then
\begin{equation}\label{graspe}\dirac(\delta_x)=\sum_{y:y\to x}\omega_{y\to x}- \sum_{y: x\to y}\omega_{x\to y},\quad \dirac(\omega_{x\to y})= \lambda_{y\to x}\delta_y-\lambda_{x\to y}\delta_x\end{equation}
and is no longer sensitive to the connection used, i.e. just depends on the metric.

Next, for a (pre)Hilbert space inner product we need an integration on $A=\C(X)$ of the form
\[ \int a=\sum_x \mu_x a(x),\quad \mu_x>0\]
so that $\int a^*a>0$ for nonzero $a\in \C(X)$. Similarly 
\[ \int(\eta^*,\eta)=\sum_{x\to y}(-\mu_y\lambda_{y\to x}) |\eta^{x,y}|^2\] for $\eta=\sum_{x\to y}\eta^{x,y}\omega_{x\to y}$. Hence this is positive definite iff 
\begin{equation}\label{lamneg} \lambda_{x\to y}<0\end{equation}
for all $x\to y$.  Moreover, the condition for $\int$ to be divergence-compatible is 
\begin{equation}\label{grariemu} \mu_y\lambda_{y\to x}=\mu_x\lambda_{x\to y}\end{equation}
for all $x\to y$. In this case,  extended trace  conditions on $\int$ in Theorem~\ref{thmspe} are empty on functions and on 1-forms just repeats (\ref{grariemu}) already imposed.  The Hilbert space structure is then
\begin{equation}\label{grahilb} \<a+\omega, b+\eta\>=\sum_x \mu_x \overline{a(x)}b(x)+ \sum_{x\to y}(-\lambda_{y\to x})\overline{\omega^{x,y}}\eta^{x,y}\end{equation}
for $\omega, \eta$ with arrow coefficients $\{\omega^{x,y}\}$ and $\{\eta^{x,y}\}$.  In summary, we have a geometrically realised spectral triple by our construction for a negative metric $\{\lambda_{x\to y}\}$ obeying (\ref{lamneg}), an inner connection defined by $
\{\sigma_{x,y}\}$ obeying (\ref{grasym}) and a measure $\{\mu_x\}$ obeying (\ref{grariemu}). This operator (\ref{graspe}) is essentially one of the operators proposed in \cite{Dav} and later works such as\cite{Bia}, but now quantum geometrically realised. Moreover, (\ref{graspe}) should be viewed as the partner of the graph Laplacian on functions and 1-forms in  \cite{Ma:gra}, necessarily recovered as $\dirac^2$ according to (\ref{dirsq}).

\begin{corollary}\label{bardirgra} Given a graph with negative arrow weights and a measure $\mu$ such that (\ref{grariemu}) holds, $\sigma=\id$ obeys the remaining condition (\ref{grasym}) and geometrically realises the Dirac operator  (\ref{graspe}). 
\end{corollary}
\proof We just take the  canonical   `bare' connection $\sigma=\id$ and hence in particular $\sigma_{x,x}=\id$ (the relevant identity metrix) as needed for (\ref{grasym}). Thus all conditions are met, and we can always geometrically realise (\ref{graspe}) in this case. \endproof

 For example, if the metric is edge-symmetric then $\mu_x=\mu_y$ i.e. a constant on each connected component of the graph solves (\ref{grariemu}).  The connection $\nabla,\sigma$ is only needed for the geometric realisation and does not directly enter (\ref{graspe}) itself. The corollary says that we always take the bare one, which does the job  but typically without the optional property (\ref{covla}) of $\la$ covariant for a full geometric realisation. 

\begin{example} \rm For a non-edge symmetric example one can take the case of the $A_n$ graph $\bullet$-$\bullet$-$\cdots$-$\bullet$ with $n$ nodes numbered  $i=1,\cdots,n$. The natural metric that admits a QRG is\cite{ArgMa3} has $\lambda_{i\to i+1}=1/h_i$ free (we take these now to be negative) and $\lambda_{ i+1\to i}={1\over \phi_i h_i}$ for certain `direction coefficients'
\[ \phi_i={(i+1)_q\over (i)_q},\quad q=e^{\imath\pi\over n+1},\quad  (i)_q={q^i-q^{-i}\over q-q^{-1}}\]
in terms of symmetric $q$-integers. We can then just take
\[ \mu_i=(i)_q \]
to solve (\ref{grariemu}) and note that these are all positive as required. We take the bare connection, with Dirac operator given by (\ref{graspe}). On the other hand, we cannot take the QLC in \cite{ArgMa3} for this metric as it does not obey (\ref{grasym}). The same applies for the QRG of $\N$ regarded as the limit $n\to\infty$. 
\end{example}

 For a fully geometric example where $\nabla$ can be taken as a QLC so that $\doublenabla(\la)=0$, we consider the $n$-gon graph with vertices identified with $\Z_n$. 
 
 \begin{example}\rm For the $n$-gon with nodes $i=0,\cdots n-1$, the natural QLC  for any edge-symmetric metric coefficients $\lambda_{i\to i+1}=\lambda_{i+1\to i}$ is reviewed in \cite{ArgMa4}   (and is unique for $n>4$) in terms of left invariant 1-forms $e^\pm:=\sum_{i}\omega_{i\to i\pm 1}$. Using $\omega_{i\to i\pm 1}=\delta_ie^\pm$ and $\lambda_{i\to i+1}=1/a(i)$ in terms of the metric coefficient function $a$ in \cite{ArgMa4}, one can find explicitly
\[ \nabla\omega_{i\to i\pm 1}=\omega_{i\pm1\to i\to i\pm 1}-\omega_{i\to i\pm 1\to i}-\omega_{i\mp 1\to i\to i\pm 1}-\rho_\pm(i)\omega_{i\to i\pm 1\to i\pm 2}\]
\[ \sigma(\omega_{i\to i\pm1\to i\pm2})=\rho_\pm(i) \omega_{i\to i\pm1\to i\pm2},
\quad \sigma(\omega_{i\to i\pm1\to i})=\omega_{i\to i\mp1\to i}, \]
where 
\[ \rho_\pm(i)={\lambda_{i\to i\pm 1}\over\lambda_{i\pm 1\to i\pm 2}}.\]
The spaces $\Omega^2_{i,i\pm 2}$ are each 1-dimensional (stepping in the same direction) and $\sigma$ acts by $\rho_\pm(i)$. The spaces $\Omega^2_{i,i}$ are each 2-dimensional (stepping back and forth in either direction) and $\sigma$ acts by swapping the direction. The general geometric Dirac operator from (\ref{gradir}) for this connection is then
\[\dirac(\delta_i)=\omega_{i-1\to i}+\omega_{i+1\to i}-\omega_{i\to i+1}-\omega_{i\to i-1},\quad \dirac(\omega_{i\to i\pm 1})=\lambda_{i\pm 1\to i}\delta_{i\pm 1}-\lambda_{i\to i\mp 1}\delta_i.\]
We also solve (\ref{grariemu}) with measure $\mu_i=1$ (say) since the metric is edge-symmetric.  However, for the geometric spectral triple construction, the condition (\ref{grasym}) holds only when $\lambda_{i\to i+1}=\lambda_{i\to i-1}$, which forces us to the regular polygon case where $\lambda_{i\to i+1}=\lambda$ independently of $i$. There are $n+2$ zero eigenvalues and the remaining eigenvalues go as $\pm\imath \sqrt{-\lambda}$ up to real constants that depend on $n$. Otherwise,  the failure of (\ref{grasym}) means that the second part of  (\ref{Jeqns}) in the construction of $\CJ$ in Lemma~\ref{lemJ} does not hold for the QLC choice of $\nabla$. In the non-constant case we can still take the bare connection as per the  corollary, with $\dirac$ according to (\ref{graspe}). \end{example}

These examples illustrate that the required $\nabla$ for our construction of $\nabla_S$ might not be a QLC for the metric. Indeed, we have the freedom on our construction to choose $\nabla$ differently if we are not aiming for the covariance of $\la$. Spectral triples in general on graphs have been analysed in \cite{Kra}.

\subsection{$2\times 2$ matrices $M_2(\C)$ case}\label{setmat}

We take $A=M_2(\C)$ with its standard 2D inner $\Omega^1$ with central basis $s,t$, $\extd=[\theta,\ \}$ for 
\[ \theta=E_{12}s+E_{21}t\]
 and wedge product relations $s\wedge t=t\wedge s, s^2=t^2=0$ (that $s,t$ commute is here forced by the structure of $\Omega^1$). The $*$ structure is hermitian conjugation in $A$ and $s^*=-t$. We take the two `standard' quantum metrics for which QLCs were studied in \cite{BegMa}, namely
\[ {\rm Case\ (i):}\quad g=\lambda^{-1}(s\tens t-t\tens s),\quad {\rm Case\ (ii):}\quad g=\lambda^{-1}(s\tens s+t\tens t)\]
where we have inserted a fixed real normalisation factor $\lambda$ for convenience later. The principal component of the moduli space of $*$-preserving QLCs in the first case is 4-dimensional in Case (i) \cite[Example~8.13]{BegMa} and  3-dimensional in Case (ii) \cite[Exercise~8.3]{BegMa}. 
 
 \begin{lemma} Among the principal components of the moduli spaces of $*$-preserving QLCs, the $\sigma$-symmetry condition (\ref{sig-sym}) holds in Case (i) for the unique point
 \[ \nabla s=2E_{12}s\tens s+2 E_{21} t\tens s,
 \quad \nabla t=2E_{12}s\tens t+2 E_{21} t\tens t\]
and in Case (ii) for the 1-parameter line
 \[ \nabla s=2 E_{21}t\tens s+ \rho E_{21}(t\tens t- s\tens s) - \rho E_{12}(s\tens t-t\tens s)\]
 \[ \nabla t=2 E_{12}s\tens t+ \rho E_{12}(t\tens t- s\tens s) - \rho E_{21}(s\tens t-t\tens s)\]
where $\rho\in \imath \R$ is imaginary.
 \end{lemma}
 \proof The equations for $\nabla$ and (\ref{sig-sym}) are not sensitive to the normalisation of the metric, so we omit the normalisation factor $\lambda$ for simplicity.  In Case (i) we have $(s,s)=0$ while $(\ ,\ )\sigma(s\tens s)=-2\alpha$ using the parameterization of QLCs in \cite[Example~8.13]{BegMa}, so $\alpha=0$ there. Similarly, $(t,t)=0$ forces $\beta=0$ there. Next $(s,t)=-1$ while  $(\ ,\ )\sigma(s\tens t)=-1-2\nu$ forces $\nu=0$ there, while $(t,s)=1$ similarly forces $\mu=0$ there. This leaves the base connection $\nabla=2\theta\tens(\ )$ on the basis, corresponding to $\sigma=-{\rm flip}$. In Case (ii) we have $(s,s)=1$ and $(\ ,\ )\sigma(s\tens s)=1-\mu-\nu$ using the parameterization of QLCs in \cite[Exercise~8.3]{BegMa}, so that $\mu+\nu=0$ there. Then $(s,t)=0$ while $(\ ,\ )\sigma(s\tens t)={\mu\over\rho}(\mu+\nu-2)$ which (for generic $\rho$ as assumed here) forces $\mu=0$. This takes us to the 1-parameter case stated, which applies for all imaginary $\rho$ including 0. This 1-parameter case was studied recently in \cite{LirMa3,BegMa:cur} and has
 \[ \sigma(s\tens s)=s\tens s+\rho(s\tens t-t\tens s),\quad \sigma(t\tens t)=t\tens t+\rho(s\tens t-t\tens s),\]
\[\sigma(s\tens t)=-t\tens s+\rho(s\tens s-t\tens t),\quad \sigma(t\tens s)=-s\tens t+\rho(s\tens s-t\tens t).\]
and vanishing Ricci curvature at $\rho=\pm\imath$ (and vanishing Ricci scalar for all $\rho$). \endproof

We proceed with these QLCs so that our construction applies and we have a geometrically realised spectral triple at the local tensorial level. In both cases $\CS$ is 3-dimensional over the algebra, with basis $1,s,t$.  Writing $\phi=\phi_0+\phi_s s+\phi_t t$, the Clifford action on  the column vector of coefficients in $M_2(\C)$ is  easily seen to be given by the matrices
\[ {\rm Case\ (i):}\quad C^s=s\la= \begin{pmatrix} 0 & 0 & -\lambda\\ 1 & 0 & 0\\ 0 &0 &0  \end{pmatrix},\quad C^t=t\la= \begin{pmatrix}   0 & 
\lambda & 0\\ 0 & 0 & 0\\ 1 & 0 & 0\end{pmatrix}\]
\[ {\rm Case\ (ii):}\quad C^s=s\la= \begin{pmatrix} 0  & \lambda & 0\\ 1 & 0 & 0\\ 0 &0 &0  \end{pmatrix},\quad C^t=t\la= \begin{pmatrix}   0 & 0& 
\lambda \\ 0 & 0 & 0\\ 1 & 0 & 0\end{pmatrix}\]
using the form of the metric in the two cases. Meanwhile, the divergence on the basis 1-forms is   
\[ {\rm Case\ (i):}\quad \delta s=2\lambda E_{21},\quad \delta t=-2\lambda E_{12},\quad {\rm Case\ (ii):}\quad \delta s=0=\delta t\]
so that in general 
\[ \delta(\phi_s s +\phi_t t)=(\extd \phi_s, s)+(\extd \phi_t, t)+\phi_s\delta s+\phi_t\delta t=\lambda\begin{cases} \{E_{21},\phi_s\}-\{E_{12},\phi_t\} & {\rm Case\ (i)}\\ [E_{12},\phi_s]+[E_{21},\phi_t] &{\rm Case\ (ii)}.\end{cases}    \]
The Dirac operator $\dirac=\extd+\delta$ is therefore
\begin{equation}\label{dirmat} \dirac\phi=[E_{12},\phi_0]s+[E_{21},\phi_0]t+\lambda\begin{cases}\{E_{21},\phi_s\}-\{E_{12},\phi_t\} &{\rm Case\ (i)}\\
[E_{12},\phi_s]+[ E_{21},\phi_t]  &{\rm Case\ (ii)}.\end{cases}\end{equation}
Note that this necessarily obeys 
\[ [\dirac, a]=\la(\extd a\tens \id+ a\nabla_S)- a\la\nabla_S=[E_{12},a] C^s+  [E_{21},a] C^t \]
 for all $a\in M_2(\C)$, as operators on $\CS$. Moreover, in Case (ii) one has 
 \[ \dirac=C^s[E_{12},\ ]+C^t[E_{21},\ ]\]
  precisely because $\delta=0$ on the basis. This has an `inner' form like some of the examples in \cite{Bar}. By construction $\CJ$ and $\gamma$ are likewise given explicitly by
\[ \CJ\phi=\phi_0^\dagger - \phi_t^\dagger s-\phi_s^\dagger t,\quad \gamma\phi=\phi_0- \phi_s s-\phi_t t.\]
This describes the local tensorial part of the construction.

 For the Hilbert space, we take $\int={1\over 2}{\rm Tr}$ as a normalised trace functional. On the other hand, \[  \int( \phi_s s+\phi_tt, \psi_s s+ \psi_t t)=\lambda{{\rm Tr}\over 2}\begin{cases}\phi_t \psi_s-\phi_s \psi_t  & {\rm Case\ (i)}\\
\phi_s\psi_s+ \phi_t\psi_t &{\rm Case\ (ii)},\end{cases}\]
which does not obey the condition extended trace condition in Theorem~\ref{thmspe} in Case (i) but does in Case (ii).  This means that the former case does not have $\CJ$ an antilinear isometry (but $\dirac$ and $\gamma$ are not affected). In terms of coefficients $
\phi=\phi_0+\phi_s s+\phi_t t\in \CS$ and similarly for $\psi\in \CS$, the inner product is
\[ \<\phi,\psi\>={{\rm Tr}\over 2}\begin{cases}\phi_0^\dagger\psi_0 -\lambda \phi_s^\dagger \psi_s +\lambda\phi_t^\dagger \psi_t& {\rm Case\ (i)}\\
\phi_0^\dagger\psi_0-\lambda\phi_t^\dagger \psi_s- \lambda\phi_s^\dagger\psi_t &{\rm Case\ (ii)}.\end{cases}\]
 In Case (i), we do not have the extended-positive-definite condition but in Case (ii) we do, provided we use $\lambda<0$. For then, choosing a self-adjoint basis and associated components
\[ s^1={s+t\over\imath},\quad s^2={s-t\over 2},\quad \phi_1={\imath\over 2}(\phi_s+\phi_t),\quad \phi_2={1\over 2}(\phi_s-\phi_t)\]
we have in Case (ii) that
\[ \<\phi,\psi\>={{\rm Tr}\over 2}\left(\phi_0^\dagger\psi_0-2\lambda(\phi_1^\dagger\psi_1+ \phi_2^\dagger\psi_2)\right),\quad\lambda<0\]
as the positive definite inner product on $\CS$. Thus, at least in Case (ii), we have a full spectral triple fully realised from QRG, including $\doublenabla(\la)=0$. 

Had we used the `bare dirac operator' with inner connection defined by $\sigma=\id$, we would have
\begin{align*} \dirac \phi&=\extd \phi_0+(\theta,\phi_s s+\phi_t t)-(\phi_s s+\phi_t t,\theta)\\ &=[E_{12},\phi_0]s+[E_{21},\phi_0]t+(\theta,s)\phi_s-\phi_s(s,\theta)+(\theta,t)\phi_t- \phi_t(t,\theta)\end{align*}
which on evaluation against $\theta=E_{12}s+E_{21}t$ in the two cases turns out to give exactly the same result as (\ref{dirmat}). Thus, as with graphs, the Dirac operator itself is not very sensitive to the connection and taking the `bare' choice gives the same result as the geometric choice based on a QLC but without covariance of the underlying $\la$ coming from metric compatibility. Indeed, the `bare' connection has no classical analogue since $\sigma=\id$ is not possible classically. 

The above results have some similarities to, but are strictly different from, the Dirac operators found on \cite{LirMa3} in $M_2(\C)$ with $\CS$ 2-dimensional. There, it was the Case (i) quantum metric which gave a spectral triple while Case (ii) resulted in non-hermitian $\dirac$. In the present paper with $\CS$ now 3-dimensional, it is Case (ii) which gave a full spectral triple while the construction for Case (i) does not have positive definite inner product. Note that our construction is very specific and we do not exclude the possibility of many other interesting spectral triples with $\CS$ 3-dimensional in either metric case. These would require a search along the lines of \cite{LirMa3} but now with bigger matrices.  We also note several other interesting spectral triples on matrix algebras coming from other contexts, such as \cite{Bar,BG} among recent work.

\subsection{Fuzzy sphere case}\label{setfuz}

The fuzzy sphere is $U(su_2)$ modulo a constant value of the quadratic Casimir, i.e., the standard quantisation of a coadjoint orbit in $su_2^*$. Explicitly, we take generators $x^i$ for $i=1,2,3$ with commutation relations and a differential calculus\cite[Example~1.46]{BegMa}
\[ [x^i,x^j]=2\imath\lambda_P \eps_{ijk}x^k,\quad \sum (x^i)^2=1-\lambda_P^2,\quad [x^i,s^j]=0,\quad \extd x^i=\eps_{ijk}x^js^k   \]
where $s^i$ are a central basis with $(s^i)^*=s^i$ and $\lambda_P$ is a dimensionless positive parameter. We sum over repeated indices. In the exterior algebra, the $s^i$ anticommute, but $\extd s^i$ is not zero. When $\lambda_P={1\over n}$ for $n$ a natural number, the standard $n$-dimensional representation of $su_2$ descends to one of the fuzzy sphere. (This leads to `reduced fuzzy spheres' isomorphic to $M_n(\C)$ if we quotient out by the kernel of this representation.) We define partial derivatives by 
\[ \extd a=(\del_i a)s^i,\quad \del_i={1\over 2\imath\lambda_P} [x_i,\ ]\]
for all $a\in A$, which classically are the killing vector fields for the action of orbital angular momentum. 

A generalised quantum metric has the form
\[ g=g_{ij}s^i\tens s^j,\quad (s^i,s^j)=g^{ij}\]
where $g_{ij}\in \C$ defines a hermitian matrix and $g^{ij}$ defines the inverse matrix. The simplest case is the `round metric' with $g^{ij}=\delta^{ij}$ on our unit sphere.  

Next, $\CS=A\oplus\Omega^1$ is 4-dimensional over the algebra and we write $\phi=\phi_0+\phi_i s^i$ with coefficients in $A$. The Clifford action is then
\[ s^i\la \phi=s^i\phi_0+g^{ij}\phi_j;\quad C^i=s^i\la=\begin{pmatrix} 0 & g^{i1} & g^{i2} & g^{i3}\\ \delta_{i1} & 0 & 0 &0\\  \delta_{i2} & 0 & 0 &0\\\ \delta_{i3} & 0 & 0 &0\end{pmatrix}\]
where we also exhibit the action as matrices acting on the coefficients $\phi_0,\phi_i$ as a 4-vector. Note that these do not represent a usual Clifford algebra but do have a memory of it in the form ${\rm Tr}(C^iC^j)=g^{ij}+g^{ji}$.

For a quantum metric, we impose quantum symmetry, which amounts to $g^{ij}$ symmetric as a matrix (and hence with real entries), and then there is a unique $*$-preserving QLC with constant Christoffel symbols\cite[Prop.~3.4]{LirMa}, namely
\[ \nabla s^i=-{1\over 2}g^{ij}( 2\eps_{jlm}g_{mk}+{\rm Tr}(g)\eps_{jkl}   )s^k\tens s^l\]
with $\sigma={\rm flip}$ on the $s^i$. With this QLC, the $\sigma$-symmetry condition reduces to quantum symmetry already imposed since $(\phi,\psi)=\phi_ig^{ij}\psi_j$ while $(\ ,\ )\sigma(\phi\tens\psi)=(\ ,\ )(\phi_i s^j\tens s^i \psi_j)=\phi_ig^{ji}\psi_j$ when $\phi=\phi_is^i\in \Omega^1$ and similarly for $\psi$. We can therefore proceed with our canonical construction at the local tensorial level. From $\nabla$ above, we see that
\[ \delta s^i=0\]
and hence the Dirac operator is 
\begin{equation}\label{dirfuz} \dirac\phi =(\del_i\phi_0)s^i+ (\extd \phi_i,s^i)+ \phi_i\delta s^i=(\del_i\phi_0)s^i+ g^{ij}\del_i\phi_j,\quad \dirac=C^i\del_i\end{equation}
where we also write it in terms of the `Clifford action' matrices acting on column 4-vectors for the coefficients. The conjugation and `even' structure are
\[ \CJ(\phi)=\phi_0^*+ \phi_i^* s^i,\quad \gamma(\phi)=\phi_0 - \phi_i s^i.\]

This completes the local tensorial level of the spectral triple. For the pre-Hilbert space structure we use $\int:A\to \C$ defined at the algebraic level as the rotationally invariant component under the action of $SU_2$ on $A$ (this has to be a multiple of 1 as the algebra has trivial centre), see \cite{LirMa}. This is a trace and, moreover, 
\[ \int(\phi,\psi)= \int \phi_i g^{ij}\psi_j=\int \psi_j g^{ji}\phi_i=\int(\psi,\phi)\]
as needed for the extended trace condition, since $\int$ is a trace and $g^{ij}$ is symmetric and $\C$-valued. Hence we have a (pre)spectral triple with
\[ \<\phi,\psi\>=\int\phi_0^*\psi_0 + \int \phi_i^* g^{ij}\psi_j\]
provided we take $g^{ij}$ with strictly positive eigenvalues so that this is positive definite. 

The calculus here is inner in degree 1 with $\theta={1\over 2\imath\lambda_P}x^is^i$ and hence there is a `bare' connection with $\sigma=\id$. This gives 
\[\dirac\phi=\extd \phi_0+ (\theta,\phi_is^i)-(\phi_i s^i,\theta)=(\del_i\phi_0)s^i+ {1\over 2\imath\lambda_P}(x^j g^{ji}\phi_i-\phi_i g^{ij}x^j)=(\del_i\phi_0)s^i+ g^{ji}\del_j\phi_i \]
which is again the same as for the geometric QLC. Although $A$ is no longer finite-dimensional, it is graded by angular momentum under the action of $SU_2$ and each graded component is finite-dimensional. As with the graph and matrix cases, $\dirac$ itself is not very sensitive to the connection used and we get the same using the `bare' connection, but without covariance of $\la$ coming from metric compatibility and with a bimodule connection that has no classical analogue. 

This Dirac operator is strictly different from the natural geometric Dirac operator on the fuzzy sphere with $\CS$ 2-dimensional over the algebra as found in \cite{LirMa2}, which has a more conventional form with Pauli matrices for the Clifford structure. The same applies to the rotationally-invariant fuzzy sphere spectral triples studied in \cite{And}. 

\section{Dirac operator with an external bundle and connection}\label{secE}

Finally, if we want to have a Dirac operator minimally coupled to an external gauge field then we should have a separate bimodule $E$ with bimodule connection $\nabla_E:E\to \Omega^1\tens_A E$. We still assume a quantum metric and bimodule connection $\nabla$ with braiding $\sigma$ on $\Omega^1$. In this case, there is an extended divergence
\[ \delta_E: \Omega^1\tens_A E\to E,\quad \delta_E=((\ ,\ )\tens\id)\nabla_{\Omega^1\tens E}\]
which given the form of the tensor product connection amounts to
\[ \delta_E(\omega\tens e)=((\ ,\ )\tens\id)\nabla\omega)e+ ((\ ,\ )\circ\sigma\tens\id)(\omega\tens \nabla_E e)=(\delta\omega)e+((\ ,\ )\circ\sigma\tens\id)(\omega\tens \nabla_E e) \]
for all $e\in E, \omega\in \Omega^1$. Then $\delta_E$ obeys 
\begin{align*} \delta_E(a\omega\tens e)&=a\delta_E(\omega\tens e)+(\extd a,\omega)e,\\
\delta_E(\omega\tens e a)&=(\delta_E(\omega\tens e)) a+((\ ,\ )\circ\sigma\tens \id)(\id\tens\sigma_E)(\omega\tens e\tens\extd a)
\end{align*}
generalising the Leibniz properties of $\delta$ before.  It is then natural to take $\CS=E\oplus \Omega^1\tens_A E$ with 
\[ \nabla_S|_E=\nabla_E, \quad \nabla_S|_{\Omega^1\tens E}=\nabla_{\Omega^1\tens E}\]
and define the Clifford action on $\CS$ by
\[ \eta\la (e+\omega\tens f)=\eta \tens e+ (\eta,\omega)f\]
for $e,f\in E$ and $\omega,\eta\in \Omega^1$, so that 
\[ \dirac_E=\la\circ\nabla_S=\nabla_E+ \delta_E\]
as a generalisation of our previous construction (which is recovered with $E=A$ and $\nabla_E=\extd$ with the usual identifications). 
 Because of the
left Leibniz property of $\delta_E$ and $\nabla_E$, we still have
\begin{equation}\label{DEcom} [\dirac_E, a](e+\omega\tens f)=\extd a\tens e + (\extd a,\omega)f=\extd a\la(e+\omega\tens f).\end{equation}
We can also still define and check that
\[ \gamma_E(e+ \omega\tens f):=e-\omega\tens f,\quad \gamma_E^2=\id,\quad \dirac_E\gamma_E=-\gamma_E\dirac_E.\]
Finally, for $\CJ_E$, we follow the same format as before and require $E$ to have its own $*$-structure defined as an antilinear skew-bimodule map $*_E:E\to E$ squaring to the identity. We inherit associated tensor product antilinear skew-bimodule maps, which we will denote generically as $\dagger$, notably
\[  \dagger={\rm flip}(*_E\tens *_E):E\tens_A E\to E\tens_A E,\quad \dagger={\rm flip}(*\tens *_E):\Omega^1\tens_A E \to E\tens_A\Omega^1,\]
 etc. as needed. We can then define 
 \[ \CJ_E|_E=*_E,\quad \CJ_E|_{\Omega^1\tens E}=\sigma_E\circ\dagger,\quad \CJ_E(e+\omega\tens f)=*_E(e)+\sigma_E(*_E(f)\tens\omega^*)\]
 as a well-defined antilinear skew-bimodule map. We used the braiding $\sigma_E$ to restore the order of factors. It is natural at this point to suppose that $\nabla_E$ is $*$-compatible in the sense
\[  \nabla_E \circ *_E=\sigma_E\circ\dagger\circ\nabla_E\]
analogously to (\ref{Jeqns}), since this implies\cite{BegMa} that $(\sigma_E\dagger)^2=\id$ and hence that $\CJ_E^2=\id$. Clearly,  $\CJ_E\gamma_E=\gamma_E\CJ_E$, and $\CJ_E b \CJ_E^{-1}\phi= \phi b^*$ also still holds as before for all $\phi\in\CS$ and $b\in A$. Given (\ref{DEcom}) amd that $\la$ is a bimodule map, we see that all of the axioms (2)-(3) for the local tensorial level hold with $\eps=\eps'=\eps''=1$ except possibly $\CJ_E\dirac_E=\dirac_E\CJ_E$. 

It remains to be seen what is the best way to handle this remaining condition, but one approach is to define $\CS^R=E\oplus E\tens_A\Omega^1$ and note that $\nabla,\nabla_E$ imply {\em right-handed} bimodule connections 
\[ \nabla^R:=\sigma^{-1}\circ\nabla,\quad \nabla^R_E:=\sigma_E^{-1}\circ\nabla_E\]
with generalised braidings given by the inverses of the left handed ones\cite[Lemma~3.7]{BegMa}. These imply a tensor product right connection on $E\tens_A\Omega^1$ and we have a parallel right handed codifferential
\[ \delta_E^R:E\tens_A\Omega^1\to E,\quad  \delta_E^R:=(\id\tens (\ ,\ ))\nabla_{E\tens\Omega^1}^R\]
\[\delta_E^R(e\tens \omega)=e\delta^R\omega+(\id\tens(\ ,\ )\circ\sigma^{-1})(\nabla_E^Re\tens\omega)\]
where $\delta^R\omega:=(\ ,\ )\nabla^R\omega$ is the right handed divergence on $\Omega^1$. We then define 
\[ \dirac_E^R(e+f\tens \omega):=\nabla_E^R e+ \delta^R(f\tens\omega)\]
 as the natural right-handed Dirac operator on $\CS^R$. This can also be cast as $\ra\circ\nabla^R_{S^R}$ for the parallel right connection restricting to the ones on $E$ and $E\tens\Omega^1$, and a right-handed Clifford action. To transfer this right-handed Dirac operator back to $\CS$ we use $\sigma_E$ which we extended as the identity map on the $E$ component of $\CS$ to give a bimodule map $\sigma_E:\CS^R\to \CS$. We  then define
 \[ \tilde\dirac_E:=\sigma_E\dirac^R_E\sigma_E^{-1}\]
as the right handed $\dirac_E^R$ viewed on $\CS$.

\begin{proposition}\label{propJEDE} If $\nabla,\nabla_E$ are $*$-preserving  and $(\ ,\ )$ obeys the metric reality condition then $\CJ_E\dirac_E= \tilde \dirac_E \CJ_E$.\end{proposition}
\proof We write $*$ also for $*_E$, and $\nabla_E e=\nabla_E^1e\tens\nabla_E^2e$ with a similar notation for $\nabla_E^R$. Then under our assumptions, we compute
\begin{align*} *\delta_E(\omega\tens e)&=e^*(\delta\omega)^*+(\nabla_E^2e)^* ((\ ,\ )\sigma(\omega\tens\nabla_E^1e))^*\\
&=e^*\delta^R(\omega^*) +(\nabla_E^2e)^*(\ ,\ )\dagger\sigma(\omega\tens\nabla_E^1e)\\
&=e^*\delta^R(\omega^*)+(\nabla_E^2e)^*(\ ,\ )\sigma^{-1}((\nabla_E^1e)^*\tens\omega^*)\\
&=e^*\delta^R(\omega^*)+ (\nabla^R_E{}^1(e^*))(\ ,\ )\sigma^{-1}(\nabla^R_E{}^2(e^*)\tens\omega^*)\\
&= \delta^R_E(e^*\tens\omega^*)
\end{align*}
where for the second equality, we used reality of the metric and $*\delta\omega=(\ ,\ )\dagger\nabla\omega=(\ ,\ )\nabla^R(\omega^*)=\delta^R(\omega^*)$ upon writing $\nabla$ being $*$-preserving in the form $\dagger\nabla=\nabla^R *$. We  then use  $\dagger\sigma=\sigma^{-1}\dagger$ for the third equality and  $\dagger\nabla_E=\nabla_E^R*$ for the fourth equality, and then recognise the answer. In this case
\begin{align*} \CJ_E\dirac_E(e+\omega\tens f)&=\CJ_E(\nabla_E e+\delta_E(\omega\tens f))=\sigma_E\dagger\nabla_E e+\delta_E^R(f^*\tens\omega^*)+\\
&=\nabla_E (e^*)+ \delta_E^R(f^*\tens\omega^*)=\nabla_E(e^*)+\dirac_E^R|_{E\tens\Omega^1}\dagger(\omega\tens f)\\
&=\nabla_E\CJ_E e+\dirac_E^R|_{E\tens\Omega^1}\sigma_E^{-1}\CJ_E(\omega\tens f)=\tilde \dirac_E\CJ_E(e+\omega\tens f)\end{align*}
using our first result and that $\nabla_E$ is $*$-preserving, the definitions of $\dirac^R_E$ and $\CJ_E$ and on noting that
\[ \tilde\dirac_E|_E=\sigma_E\dirac_E^R|_E=\dirac_E|_E=\nabla_E,\quad \tilde\dirac_E|_{\Omega^1\tens E}=\dirac^R_E|_{E\tens\Omega^1}\sigma^{-1}_E=\delta^R_E\sigma^{-1}_E.\]
since $\sigma_E$ was extended to act as the identity on $E$. \endproof

This reduces in the case without $E$ in the setting of Section~\ref{secdir} to $\CJ\dirac=\tilde\dirac\CJ$ where $\tilde\dirac=\dirac^R=\extd+\delta^R$. Hence $\CJ$ commutes with $\dirac$ if $\delta^R=\delta$, which happens when the quantum metric is $\sigma$-symmetric as needed in that section to get a spectral triple. Also note that Proposition~\ref{propJEDE} can be phrased more symmetrically if we abandon $\CJ_E$ and work instead with $\dagger$ extended as $*_E$ to all of $\CS$,
\[ \dagger:\CS\to \CS^R,\quad \dagger (e+\omega\tens f)=*_E(e) + *_E (f)\tens \omega^*\]
which is more similar to the map $\CJ$ in Section~\ref{secdir}. Then $\CJ_E=\sigma_E\circ\dagger$ when $\sigma_E$ is also viewed on $\CS$ by the identity on $E$, and in these terms Proposition~\ref{propJEDE} simply says that
\begin{equation}\label{dirLR} \dagger \circ\dirac_E= \dirac_E^R\circ \dagger.\end{equation}
In summary, we still obtain a spectral triple at the local tensorial level, except we saw that 
the $\CJ_E, \dirac_E$ relation needs to be modified as in Proposition~\ref{propJEDE} or the more symmetric (\ref{dirLR}).  

Finally, we also suppose an inner product $\(\ ,\ )_E: E\tens_A E\to A$ on $E$ with a reality property $*\circ(\ ,\ )_E=(\ ,\ )_E\circ\dagger$.  It is then natural to set
\[ \<e,f\>=\int (*_E(e),f)_E,\quad \<\omega\tens e,\eta\tens f\>=\int(*_E(e),(\omega^*,\eta)f)_E\]
and one can check that this is well-defined on $\CS$, that the action of $a^*$ is adjoint to the action of $a$ as required, and that $\<\ ,\ \>$ is conjugate-symmetric for any $*$-preserving linear map $\int: A\to \C$. We can add here positivity requirements, which have a  similar flavour to our extended-positive definite assumptions before, so as to have a pre-Hilbert space.  It is easy to check that $\gamma$ is then hermitian, but in general  $\CJ_E$ fails to an antilinear isometry at least if we assume the natural extended-trace properties for $\int$. Motivated by our experience with $\dirac_E$, we introduce $\<\ ,\ \>^R$ on $\CS^R$ by
\[ \<e,f\>^R:=\int (*_E(e),f)_E=\<e,f\>,\quad \<e\tens\omega, f\tens \eta\>^R:=\int (\omega^*,(*_E(e),f)_E\eta)\]
which is unchanged on the $E$ component but suitably flipped on the $E\tens_A\Omega^1$ component. We then use the extended  $\sigma_E^{-1}$ to transfer this back to $\CS$ as 
\[ \tilde\<\phi,\psi\tilde\>:=\<\sigma_E^{-1}\phi,\sigma_E^{-1}\psi\>^R\]
for all $\phi,\psi\in S$. 

\begin{lemma} Under our assumptions above, if $\int$ obeys the two extended-trace conditions 
\[ \int (e,f)_E=\int (f,e)_E,\quad \int(\omega,\eta)=\int(\eta,\omega)\]
then $\tilde\<\CJ_E\phi,\CJ_E\psi\tilde\>=\<\psi,\phi\>$ or, equivalently, $\<\phi^\dagger,\psi^\dagger\>^R=\<\psi,\phi\>$, for all $\phi,\psi\in \CS$.\end{lemma}
\proof Between the $E$ components and writing $*$ for $*_E$, it is immediate that $\<e^*,f^*\>^R=\int (e,f^*)_E=\int (f^*,e)_E=\<f,e\>$. We also need
\begin{align*}\<\dagger(\omega\tens e),\dagger(\eta\tens f)\>^R&=\<e^*\tens \omega^*,f^*\tens \eta^*\>=\int (\omega,(e,f^*)_E\eta^*)=\int (\omega(e,f^*)_E,\eta^*)\\
&=\int (\eta^*,\omega (e,f^*)_E)=\int (\eta^*,\omega )(e,f^*)_E=\int ((\eta^*,\omega )e,f^*)_E\\
&=\int (f^*,(\eta^*,\omega )e)_E =\<\eta\tens f,\omega\tens e\>\end{align*}
using the two extended-trace properties for the 4th and 7th equalities, and that the inner products are bimodule maps and defined on the tensor product over $A$. 
\endproof
There is a similar issue with $\dirac_E$ not being hermitian which we will address in more detail elsewhere. In the remainder, we  consider what the above construction amounts to in the simplest case, namely $E=A$ as a bimodule by left and right multiplication.  

\begin{proposition}\label{nablaalpha} Let $(\Omega^1,\extd)$ be inner and equipped with a quantum metric and a given bimodule connection $\nabla$ with generalised braiding $\sigma$. 

(1) Bimodule connections  $\nabla_A$ on $E=A$ are classified by $\alpha\in \Omega^1$ of the form  $\alpha=\zeta(\theta)+\alpha_0$ for some bimodule map $\zeta:\Omega^1\to \Omega^1$ and some central  $\alpha_0\in \Omega^1$. Then 
\[ \nabla_A a= (\extd a+ a\alpha)\tens 1,\quad \sigma_A(a\tens \extd b)=a(\extd b+[b,\alpha])\tens 1 \]
for all $a,b\in A$. Here $\sigma_A$ extends to $\sigma_A(a\tens\omega)=a(\id-\zeta)\omega\tens 1$ for all $\omega\in \Omega^1$.

(2)  $\Omega^1\tens_AA=\Omega^1$ via the standard identification acquires a tensor product bimodule connection
\[ \nabla_\alpha \omega=\nabla\omega+\sigma(\omega\tens \alpha),\quad \sigma_\alpha(\omega\tens\eta)=\sigma(\omega\tens (\id-\zeta)\eta)\]
for all $\omega,\eta\in \Omega^1$.
\end{proposition}
\proof (1) Since $\Omega^1$ is inner, a bimodule connection by \cite{MaSim} takes the form 
\[\nabla_A a=\theta \tens a-\sigma_A(a\tens\theta)+ \alpha_A(a)\]
for freely chosen bimodule maps $\sigma_A$ and $\alpha_A$. In our case a bimodule map $\sigma_A(a\tens \omega)=a\sigma_A(1\tens \omega)=\sigma_A(1\tens a\omega)$ needs to have the form $\sigma_A(a\tens \omega)=a(\id-\zeta)\omega\tens 1$ for some freely chosen bimodule map $\zeta:\Omega^1\to \Omega^1$ (we just need $\id-\zeta$ a bimodule map but have chosen to split it this way). Similarly $\alpha_A(a)=a\alpha_A(1)=\alpha_A(1)a$ needs $\alpha_A(a)=a\alpha_0\tens 1$ for a freely chosen central element $\alpha_0\in \Omega^1$.  The resulting connection depends only on $\alpha=\zeta(\theta)+\alpha_0$ since
\[ \nabla_A a=(\theta a-a\theta+ a\zeta(\theta)+a\alpha_0)\tens 1=(\extd a+a\alpha)\tens 1,\]
\[ \sigma_A(a\tens\extd b)=a(\extd b\tens 1-\zeta([\theta,b]))\tens 1=a(\extd b-[\zeta(\theta),b])=a(\extd b+ [b,\alpha])\tens 1.\] 
Here, $\nabla_A$ depends only on $\alpha$ and hence so does $\sigma_A$, but the existence of $\zeta,\alpha_0$ is needed for $\sigma_A$ to be well-defined, i.e. to have a bimodule connection. 

(2) Once we have $\nabla_A$ as a bimodule connection, we can compute the tensor product bimodule connection \cite{BegMa}, in our case on $\Omega^1\tens_A A$.  Then 
\[\nabla_{\Omega^1\tens A}(\omega\tens a)=\nabla\omega\tens a+ (\sigma\tens\id)(\omega\tens\nabla_A a)=\nabla\omega\tens a+\sigma(\omega\tens(\extd a+ a\alpha))\tens 1\]
\[ \sigma_{\Omega^1\tens A}(\omega\tens a\tens\eta)=(\sigma\tens\id)(\omega\tens a(\id-\zeta)\eta\tens 1)=\sigma(\omega\tens a(\id-\zeta)\eta)\tens 1\]
We identify $\Omega^1\tens_AA=\Omega^1$ in the standard way, in which case it suffices to set $a=1$ to give the result stated. We denote this as $\nabla_\alpha$ since, after this identification,  the result looks like $\nabla$ modified by the 1-form $\alpha$.  \endproof 

 
With the standard identifications,  the general $\delta_E$ and $\dirac_E$ construction above just amounts to Section~\ref{secdir} with $\nabla_A,\nabla_\alpha$ in place of $\extd,\nabla$, so $\nabla_S|_A=\nabla_A$ and $\nabla_S|_{\Omega^1}=\nabla_\alpha$. The Clifford action and quantum metric are unchanged and hence
\begin{equation}\label{diracalpha} \dirac_\alpha(a +\omega)=\nabla_A a+  (\ ,\ )\nabla_\alpha\omega=\dirac(a+\omega)+a\alpha+ (\omega,\alpha)\end{equation}
as the Dirac operator modified by $\alpha$ according to our construction, assuming for convenience that the quantum metric was $\sigma$-symmetric to begin with. The bundle $E$ here is trivial (a rank 1 free module) and we can think of $\alpha$ as in the role of a $U(1)$ gauge field in this expression. Note that this modification is different from what we considered before in Proposition~\ref{a0alpha}, since we are not keeping $\sigma_S$ fixed as in the discussion there.

We also take $*_E=*$ on $A$ and $(a,b)_E=ab$ as the obvious choices for $E=A$. The sesquilinear inner product on $\CS$ and $\gamma_E$ are, like the Clifford action, unchanged from original construction in Section~\ref{secdir} once we identify $\Omega^1\tens_AA=\Omega^1$. The connection $\nabla_A$ is not necessarily $*$-preserving but if we restrict $\alpha$ so that it is then $\CJ$ gets modified to
\[ \CJ_\alpha(a+\omega)=a^*+(\id-\zeta)(\omega^*)\]
after the usual identification, according to our general constructions above. Our point of view is that even if $\CJ_\alpha, \dirac_\alpha$ have modified properties making them no longer part of a usual spectral triple, they are the natural objects coming from the general theory. 

We now specialise this simplest choice of $E$ to a graph calculus on $A=\C(X)$. Then there are no central elements $\alpha_0\in \Omega^1$ and any bimodule map $\zeta:\Omega^1\to \Omega^1$ has to have the form $\zeta(\omega_{x\to y})=\alpha_{x\to y}\omega_{x\to y}$ for some coefficients $\alpha_{x\to y}$. These coefficients define a 1-form $\zeta(\theta)=\sum_{x\to y}\alpha_{x\to y}\omega_{x\to y}=\alpha$ in the same conventions as we used before for the coefficients of a 1-form in the arrow basis. Thus, the conditions on $\alpha$ in Proposition~\ref{nablaalpha} are automatic and the data for the modified connection and Dirac operator are just any 1-form $\alpha$, i.e. an assignment of coefficients to every arrow. Then
\begin{equation}\label{graalpha} \nabla_{\alpha}(\omega_{x\to y})=\nabla\omega_{x\to y}+\sum_{y\to z}\alpha_{y\to z}\sum_{x\to w\to z}\sigma_{x,z}{}^w{}_y \omega_{x\to w\to z}\end{equation}
is the `minimally coupled' covariant derivative. This is not necessarily $*$-preserving even if $\nabla$ is and $\alpha^*=-\alpha$ (say). The modified Dirac operator in the case of the metric initially $\sigma$-symmetric is  then
\begin{equation}\label{dirgraalph} \dirac_{\alpha}(\delta_x)=\dirac(\delta_x)+\sum_{y:x\to y}\alpha_{x\to y}\omega_{x
\to y},\quad \dirac_\alpha(\omega_{x\to y})=\dirac(\omega_{x\to y})+\lambda_{x\to y}\alpha_{y\to x}\delta_x\end{equation}
where $\dirac$ is the original graph Dirac operator, for example (\ref{graspe}) for a suitable $\nabla$. We  consider this as the natural Dirac operator, of the type we had before, but now minimally coupled to an external 1-form $\alpha$. The associated charge conjugation map is 
\begin{equation}\label{Jgraalph} \CJ_\alpha(\delta_x)=\delta_x,\quad \CJ_\alpha(\omega_{x\to y})=(\alpha_{y\to x}-1)\omega_{y\to x}\end{equation}
extended antilinearly.

\section{Concluding remarks}\label{secrem}

We have associated a Dirac operator within quantum Riemannian geometry to every bimodule connection $\nabla$ with generalised braiding $\sigma$ and a $\sigma$-symmetric quantum metric $(\ ,\ )$. In particular, we can take the canonical `bare' connection on an inner $\Omega^1$ defined by $\sigma=\id$, for any quantum metric, hence giving a canonical $\dirac$ determined by the quantum metric alone. On a graph, the quantum metric consists of `square lengths' on every arrow and we essentially recover a graph Dirac operator associated to such data in \cite{Dav,Bia} and other previous works. We also illustrated the construction on $M_2(\C)$ and the fuzzy sphere, obtaining Dirac operators with respectively 3-dimensional and 4-dimensional spinor bundles over the algebra. When there is a suitable trace $\int: A\to \C$,  the general construction also obeys the pre-Hilbert space level of a spectral triple, giving good contact with the latter approach to noncommutative geometry. Finally, we showed how one can naturally `minimally couple' the construction to an external bimodule connection $\nabla_E$ and showed how this works in the simplest case of $E=A$, and in particular in the graph case where the data for the external connection just amounts to a 1-form $\alpha\in \Omega^1$, i.e. a collection $\{\alpha_{x\to y}\}$ of coefficients associated to the arrows. 

There are several  independent questions that arise from the present work. The first is that our treatment was essentially algebraic and it would be interesting to extend this to an analytic setting to handle the extension of maps to the required $L^2$ completion in the infinite-dimensional case. While this may be clear enough for the final product (constructed algebraically on a dense subalgebra) it would also be better to allow $A$ to be an operator algebra within the QRG itself. There are numerous technical issues for this, however, particularly in the handling of tensor products over $A$, putting this beyond our scope here. 

Another question, even at the algebraic level, is how to extend the construction to $\CS=\Omega$ the total differential graded or `exterior' algebra of forms on $A$. If one has a Hodge operator, or at least a codifferential $\delta$ in all degrees, then one can follow the same strategy with $\dirac=\extd+\delta$ in all degrees, and thereby aim for a noncommutative version of the classical spectral triple of this type in \cite{Con:spe}. The problem here is that quantum Riemannian geometries in the existing framework do not necessarily admit natural Hodge or higher codifferential structures in much generality. Recent work which indeed studied $\delta$ in a certain extension theory context was in \cite{Ma:rec}. Likewise, when $A$ is a quantum group with a bicovariant calculus then $\Omega$ is a super-Hopf algebra allowing one to define a kind of Hodge operator in nice cases\cite{Ma:hod}. There are also specific models, such as $q$-Minkowski space or the Bruhat graph on $S_3$ where a reasonable Hodge operator can be exhibited on the exterior algebra by ad hoc methods. The latter example is recapped in \cite[Example~1.74]{BegMa}. Note that working with an exterior algebra is very different in the graph case from working in a more usual way with a simplicial complex as in \cite{Bia}, but there could be insights from here also, and conversely applications to topological data networks. 

Finally, the construction $\dirac_E$ in Section~\ref{secE} can be taken much further. One can compute $\dirac_E$ for more complicated choices of bimodules $E$, including for example the $q$-monopole bundles on the standard $q$-spheres associated to $\C_q[SU_2]$. Likewise, on graphs, a general bimodule $E$ is just a collection ${}_xE_y$ of vector spaces, for every pair of vertices, and this can be elaborated to include a connection $\nabla_E$, effectively generalising the  1-form $\alpha$ to connections on more nontrivial vector bundles. Beyond this, one should look generally at vector bundles and connections $E,\nabla_E$ associated to quantum principal bundles and `spin connections' on them (as for the $q$-monopole). Moreover, in all these cases, as seen even for the simplest choice $E=A$, the abstract properties of $\dirac_E$ remain to be fully developed as an appropriate generalisation of standard spectral triples. The results so far suggest that the most elegant way to do this could be in a left-right symmetric manner with the role of $\CJ_E$ replaced by $\dagger: \CS\to \CS^R$ as in (\ref{dirLR}). 

It should be stressed that the classical geometry behind the type of spectral triples in the present work is very different from spinor geometry. Hence these spectral triples have a different flavour from ones that deform the classical case, such as on fuzzy spheres\cite{And,LirMa2,Bar} and noncommutative tori\cite{Bon,BG,LirMa3}. On the other hand, a topical application of even finite spectral triples (where $A$ and $\CS$ are finite-dimensional as vector spaces) is in the almost commutative case where these are tensored onto classical spacetime as an approach to particle physics \cite{Con0}. We also note \cite{DabSit} for an example of recent work in this direction. In this context, new approaches to finite examples such as on $M_2(\C)$ in the present work could be of interest.   These are some directions for further work.

\end{document}